# Numerical solution of variable order fractional differential equations.

## John T. Katsikadelis

Department of Civil Engineering
National Technical University
Athens, Greece
e-mail: jkats@central.ntua.gr

## Abstract

A method for the numerical solution of variable order (VO) fractional differential equations (FDE) is presented. The method applies to linear as well as to nonlinear VO-FDEs. The Caputo type VO fractional derivative is employed. First, an simple expression, which approximates the VO fractional derivative, is established and then a procedure based on this approximation is developed to solve VO-FDEs linear and nonlinear, both explicit and implicit. VO-FDEs with variable coefficients are also treated. The method is illustrated by solving the second order VO-FDE describing the response of the VO fractional oscillator, linear and nonlinear (Duffing). However, it can be straightforwardly extended to higher order VO-FDEs. The presented method, in addition to its effectiveness, is simple to implement and program on a computer. The obtained results validate the efficiency and accuracy of the developed method

## 1. Introduction

The formulation of physical problems reflects the mathematical tools available at the time of their development. Therefore, much of modeling of the physical world has been described by differential equations involving integer order derivatives as defined by Leibniz and Newton. However, research carried out in recent years has pointed out that fractional order derivatives provide an effective tool to reliably model many complex physical and engineering systems. This fact gave a great boost to the study of fractional differential equations governing the response of systems modeled in this way. In this context, many books have been written on Fractional Calculus as well as numerous publications on the study of physical and engineering systems via fractional derivatives, see for example in [1-5] and the references therein. It is noteworthy that though Fractional Calculus is as old as Classical Calculus, it has only recently been employed to study the response of physical and engineering systems. Historically, the first written information on fractional derivatives is in a letter written by L'Hôpital to Leibniz after he defined the integer order derivative $\partial^n y / \partial x^n$. There he asks: "*What if $n$ is a fraction, say $n = 1/2$ ?*" Leibnitz in his answer (30 September 1695) [6] gave an explanation to this question concluding: "*Ainsi il s'en suit que $d^{1:2}x$ sera égal à $x \cdot \sqrt[2]{dx:x}$*" and added prophetically "*Il y a de l'apparence qu'on tirera un jour des conséquences bien utiles de ces paradoxes, car il n'y a guerres de paradoxes sans utilité*". For 3 centuries, the fractional derivative inspired pure theoretical mathematical developments useful only for mathematicians. The integer order derivative allows giving geometrical interpretations to the proposed physical models resulting from Newton's laws [7,8]. Apparently, this made the then revolutionary concepts accessible to the contemporary scientists, who were well experienced in geometry. In a sense, the long delay to apply Fractional Calculus may be attributed to this fact. We recall that fractional derivatives lack the straightforward geometrical interpretation of their integer counterparts [9].

The fractional calculus has allowed the definition of any order fractional derivative (FD), real or imaginary. This fact enables us to consider the fractional derivative to be a function of time (explicit VO-FD) or of some other dependent variable (implicit VO-FD). Although the extension from constant-order to VO-FD may seem somewhat natural and several systems have been modeled with VO-FD (see for example [10-15]), this idea has been forwarded only very slowly [12]. In a sense, this extension may be have been prevented by the difficulty to obtain solutions to VO-FDEs. Nevertheless, VO-Calculus, besides the suitable modeling of actual structures, may model nonlinear a response in constant order differential equations as linear response in a VO-Calculus framework, with all the simplifications that arise from the use of linear operators [15]. This particular feature of VO-Calculus alone is sufficient reason to consider it as a promising alternative basis to constant order (integer or fractional) Calculus for modelling complex physical phenomena involving memory effects.

The modeling of physical systems using VO-FDs leads to VO-FDEs. Analytical solutions as in the case of constant order FDEs, are available only for few simple FDEs, therefore inefficient to obtain results for applications to real world problems. Though researchers have developed efficient numerical solution methods for ordinary and partial FDEs of constant and distributed fractional order, e.g., [16-19], the literature on numerical approximation of VO fractional derivative and numerical solution of VO-FDEs is limited, e.g. [20-24], and simple general efficient numerical methods for the solution of VO-FDEs are to author's knowledge unavailable.

In this paper, first a robust numerical expression approximating the VO-FD is established and then it is employed to develop an effective numerical method for the solution of linear and nonlinear VO FDEs, both explicit and implicit. The method applies also to VO-FDEs with variable coefficients. The method is demonstrated by solving the second order VO-FDE describing the response of the VO fractional oscillator, linear and nonlinear (Duffing). However, it can be straightforwardly extended to higher order VO-FDEs and in conjunction with the AEM (Analog Equation Method) [18] to partial VO-FDEs.

The developed numerical scheme is stable, accurate and simple to implement. There is no computational complexity in constructing the solution algorithms, while they are easy to program. The convergence and the accuracy of the method are demonstrated through well corroborated numerical examples. The obtained results validate the effectiveness of proposed method.

## 2    Approximation of the VO-FD

Without excluding other types of VO FDs, we adopt the Caputo type VO order FD of a function $u(t)$

$$D^{\alpha(t)}u(t) = \frac{1}{\Gamma[1-\alpha(t)]} \int_0^t (t-x)^{-\alpha(t)} \dot{u} dx, \qquad 0 < \alpha(t) < 1, \ t > 0 \qquad (2.1)$$

and

$$\lim_{a \to 1} D^{\alpha(t)} u(t) = \dot{u}(t), \qquad \lim_{a \to 0} D^{\alpha(t)} u(t) = u(t) - u_0 \qquad (2.2)$$

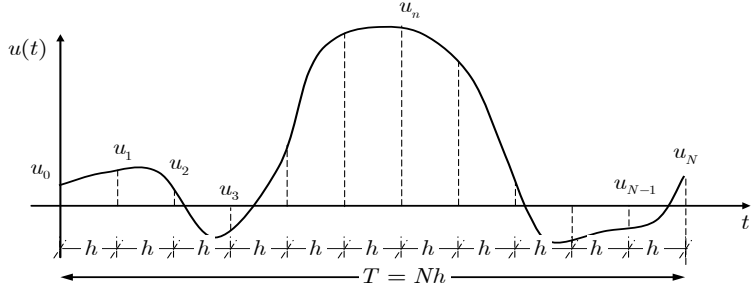

Figure 1. Discretization of the interval $[0,T]$ into $N$ equal intervals $h = T/N$

The function is considered in an interval $[0,T]$. To approximate the integrals in Eqs (1), the interval $T$ is divided into $N$ equal intervals $\Delta t = h$, $h = T/N$ (Fig. 1), in which $\dot{u}(t)$ is assumed to vary according to a certain law, e.g. constant, linear etc. In this analysis, $\dot{u}(t)$ is assumed to be constant and equal to the mean value in the interval $h$. Thus it is

$$\dot{u}_r^m = \frac{\dot{u}_{r-1} + \dot{u}_r}{2}, \qquad r = 1,2,...,N \tag{2.3}$$

Hence, the integral in Eq. (2.1) at instant $t = nh$ is written as

$$(D^{\alpha(t)}u)_n = c_1^n \dot{u}_1^m + c_2^n \dot{u}_2^m \ldots + c_n^n \dot{u}_n^m \tag{2.4}$$

where

$$\begin{aligned}
c_1^n &= \frac{1}{\Gamma[1-\alpha(nh)]} \int_0^h (nh-x)^{-\alpha(nh)} dx \\
c_2^n &= \frac{1}{\Gamma[1-\alpha(nh)]} \int_h^{2h} (nh-x)^{-\alpha(nh)} dx \\
&\ldots \\
c_r^n &= \frac{1}{\Gamma[1-\alpha(nh)]} \int_{(r-1)h}^{rh} (nh-x)^{-\alpha(nh)} dx, \quad r = 1,2,...,n \quad n = 1,2,...,N
\end{aligned} \tag{2.5}$$

**Example 1**

Compute the VO-FD of the function $u(t) = t^2$, $t \in [0,1]$, (i) $\alpha = (50t+49)/100$ (ii) $\alpha = 1 - \exp(-t)$.

(i) The exact VO-FD is

$$D^{\alpha(t)}(t^2) = \frac{1}{\Gamma[1-\alpha(t)]} \int_0^t (t-x)^{-\alpha(t)} \dot{u} dx = \frac{1}{\Gamma[1-\alpha(t)]} \frac{20000 t^{\frac{151}{100} - \frac{t}{2}}}{(50t-151)(50t-51)} \tag{1}$$

Fig. 2 shows the exact versus the approximate VO FD as computed using Eq. (2.4) and the error for $h = 0.001$. $\max(|\text{error}|) = 9.21514\text{e-}4$

(ii) The exact VO-FD is

$$D^{\alpha(t)}(t^2) = \frac{1}{\Gamma[1-\alpha(t)]} \int_0^t (t-x)^{-\alpha(t)} \dot{u} dx = \frac{1}{\Gamma[1-\alpha(t)]} \frac{2\exp(2t) t^{\exp(-t)+1}}{\exp(t)+1} \quad (2)$$

Fig. 3 shows the exact versus the approximate VO-FD and the error for $h = 0.001$,. $\max(|\text{error}|) = 4.62050\text{e-}05$

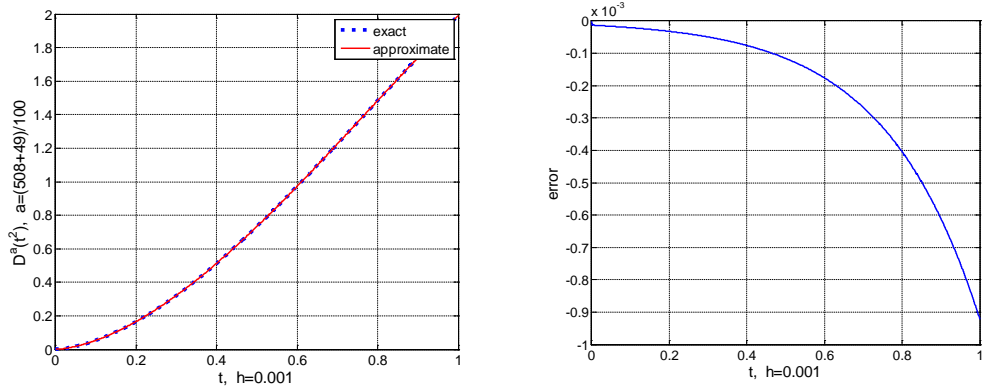

Figure 2. VO-FD in Example 1 (i).

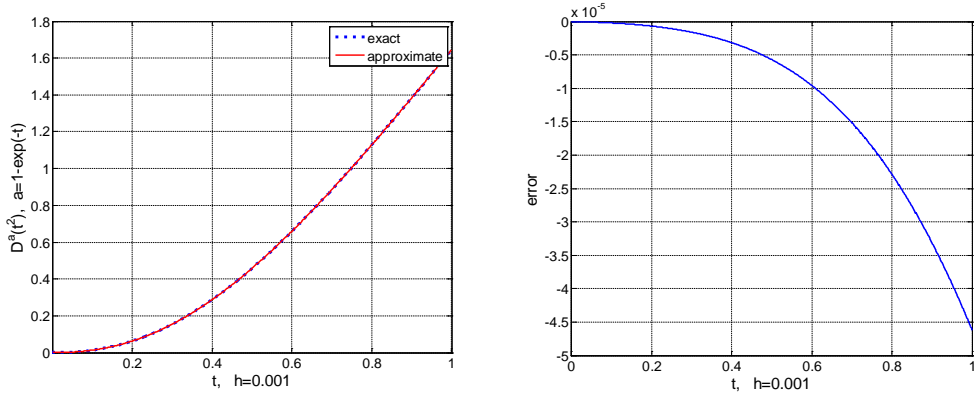

Figure 3. VO-FD in Example 1 (ii).

## 3   Linear VO-FD Equations

### 3.1. Explicit VO-DF equations

The solution procedure is illustrated by solving the second order VO-FDE describing the response of the oscillator with VO fractional damping

$$a_1 \ddot{u} + a_2 D^{a(t)} u + a_3 u = p(t), \quad u(0) = u_0, \quad \dot{u}(0) = \dot{u}_0 \quad (3.1)$$

Using the AEM as for the second order differential equation, we obtain [25]

$$u_n + \frac{c_1}{2} q_n = u_{n-1} + h\dot{u}_n - \frac{c_1}{2} q_{n-1} \quad (3.2)$$

$$\dot{u}_n - \frac{c_2}{2}q_n = \dot{u}_{n-1} + \frac{c_2}{2}q_{n-1} \tag{3.3}$$

where $q_n = \ddot{u}_n$, $c_1 = \frac{h^2}{2}$, $c_2 = h$

Applying Eq. (3.1) for $t = t_n = nh$ gives

$$a_1 \ddot{u}_n + a_2 (D^{a)}u)_n + a_3 u_n = p_n \tag{3.4}$$

Using Eq. (2.3), Eq (2.4) is written

$$\begin{aligned}(D^\alpha u)_n &= \sum_{r=1}^{n-2} c_r^n \dot{u}_r^m + c_{n-1}^n \frac{\dot{u}_{n-2} + \dot{u}_{n-1}}{2} + c_n^n \frac{\dot{u}_{n-1} + \dot{u}_n}{2} \\ &= \sum_{r=1}^{n-2} c_r^n \dot{u}_r^m + \frac{c_{n-1}^n}{2} \dot{u}_{n-2} + (\frac{c_{n-1}^n}{2} + \frac{c_n^n}{2})\dot{u}_{n-1} + \frac{c_n^n}{2} \dot{u}_n\end{aligned} \tag{3.5}$$

Substituting Eq. (3.5) in Eq. (3.4) gives

$$a_1 \ddot{u}_n + a_2 \frac{c_n^n}{2} \dot{u}_n + a_3 u_n = -a_2(\frac{c_{n-1}^n}{2} + \frac{c_n^n}{2})\dot{u}_{n-1} + g_n \tag{3.6}$$

where

$$g_n = p_n - a_2 \left[ \sum_{r=1}^{n-2} c_r^n \dot{u}_r^m + \frac{c_{n-1}^n}{2} \dot{u}_{n-2} \right], \quad c_0^1 = 0, \quad \sum_{r=1}^{n-2} c_r^n \dot{u}_r^m = 0 \text{ if } n \leq 2 \tag{3.7}$$

Eqs. (3.2), (3.3) and (3.6) are combined to a system of simultaneous equations for $q_n, \dot{u}_n, u_n$, which in matrix form are written as

$$\begin{bmatrix} a_1 & a_2 \frac{c_n^n}{2} & a_3 \\ \frac{c_1}{2} & -h & 1 \\ -\frac{c_2}{2} & 1 & 0 \end{bmatrix} \begin{Bmatrix} q_n \\ \dot{u}_n \\ u_n \end{Bmatrix} = \begin{bmatrix} 0 & -a_2(\frac{c_{n-1}^n}{2} + \frac{c_n^n}{2}) & 0 \\ -\frac{c_1}{2} & 0 & 1 \\ \frac{c_2}{2} & 1 & 0 \end{bmatrix} \begin{Bmatrix} q_{n-1} \\ \dot{u}_{n-1} \\ u_{n-1} \end{Bmatrix} + \begin{Bmatrix} g_n \\ 0 \\ 0 \end{Bmatrix} \tag{3.8}$$

Eq. (3.8) is solved successively for $n = 1, 2, \ldots$ to yield the solution $u_n$ and the derivatives $\dot{u}_n$, $q_n = \ddot{u}_n$ at instant $t = nh \leq T$. For $n = 1$, the value $q_0$ appears in the right hand side of Eq. (3.8). This value can be obtained from for (3.4) for $t = 0$. This yields

$$q_0 = (p_n - a_2(D^{a)}u)_0 - a_3 u_0)/a_1 \tag{3.9}$$

$u_0, \dot{u}_0$ are known from the initial conditions. $(D^a u)_0 = 0$ if the integrand in (2.1) is continuous, otherwise $(D^a u)_0$ can be approximated using the method developed in [26]. The problem can be also solved efficiently using the procedure for implicit VO FDEs described in Section 3.2.2.

### 3.1.1 Stability of the numerical scheme

Stability is ensured if the spectral radius $A_n$ is satisfies the condition [25]

$$\rho(A_n) \leq 1 \tag{3.10}$$

where $A_n$ is the amplification matrix

$$A_n = \begin{bmatrix} a_1 & a_2\dfrac{c_n^n}{2} & a_3 \\ \dfrac{c_1}{2} & -h & 1 \\ -\dfrac{c_2}{2} & 1 & 0 \end{bmatrix}^{-1} \begin{bmatrix} 0 & -a_2(\dfrac{c_{n-1}^n}{2}+\dfrac{c_n^n}{2}) & 0 \\ -\dfrac{c_1}{2} & 0 & 1 \\ \dfrac{c_2}{2} & 1 & 0 \end{bmatrix} \tag{3.11}$$

This condition is satisfied as it shown in Example 2.

On the base of the previously solution procedure a computer program has been written in Matlab for the solution of the initial value problem (3.1).

### Example 2

Solve the initial value problem

$$a_1\ddot{u} + a_2 D^{a(t)}u + a_3 u = p(t), \quad u_0 = 1, \; \dot{u}_0 = 10, \; a(t) = d - k\exp(-t) \tag{1}$$

Assume: $a_1 = 1, a_2 = 2\xi\omega, a_3 = \omega^2, \xi = 0.1, \omega = 5$, $p(t) = 0$ and $d, k$ as follows:

(i) $d = 0.9999, k = 10e - 10$. In this case it is $a(t) \approx 1$, hence as anticipated $D^{a(t)}u \approx \dot{u}$, and the computed solution approximates the exact solution

$$u_{ex} = \exp(-\xi\omega t)(\dfrac{\dot{u}_0 + \xi u_0}{\omega_D}\sin\omega_D t + u_0 \cos\omega_D t), \; \omega_D = \omega\sqrt{1-\xi^2} \tag{2}$$

This is shown in Fig 4.

(ii) $d = 1e - 10, k = 1e - 10$. In this case, it is $a(t) \approx 0$, hence as anticipapted $D^{a(t)}u \approx u - u_0$, and the computed solution approximates the exact solution

$$u_{ex} = \dfrac{\dot{u}_0}{\omega}\sin\omega t + u_0\cos\omega t + a_2\dfrac{u_0}{k}, \quad k = a_2 + a_3, \quad \omega = \sqrt{k/a_1} \tag{3}$$

This is shown in Fig 5.

(iii) Solve Eq. (1) for: $a = 1$; $a = 1 - \exp(-t)$; $a = 0.8$; $a = 0.8[1-\exp(-t)]$; $a = 0.5[1-\exp(-t)]$. The computed results are plotted in Fig. 6. As anticipated, the VO FD yields larger displacements than the corresponding constant order FD.

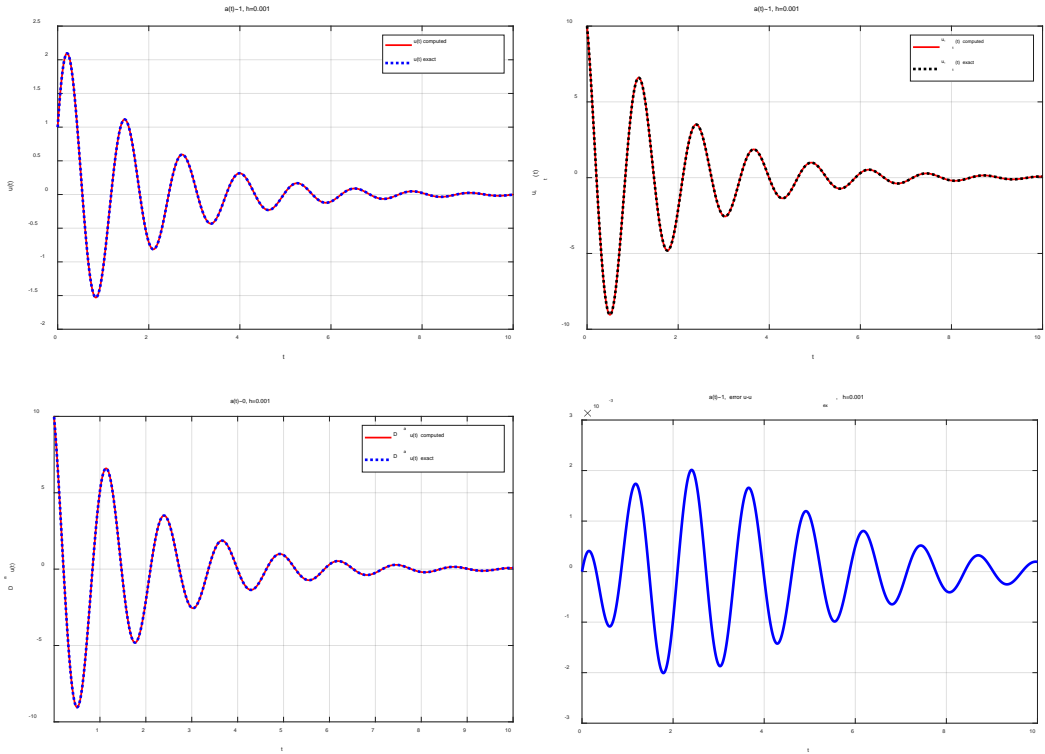

Figure 4. Results in Example 2: case (i)

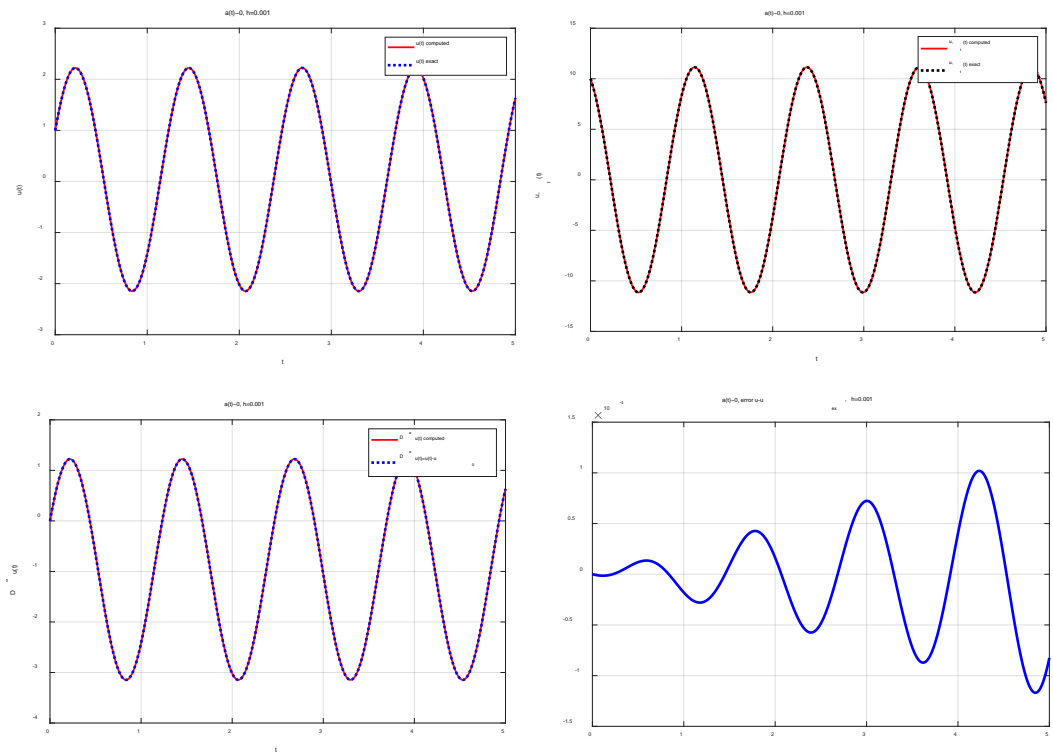

Fig. 5 Results in Example: case (ii)

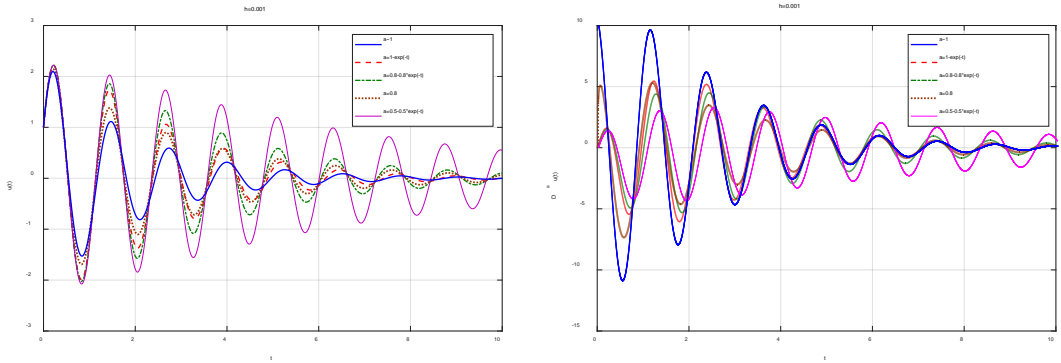

Fig. 6 Results in Example 2: case (iii)

In all studied cases the stability condition, Eq. (3.11), was satisfied. Fig. 7 shows the spectral radius for $a = 0.8[1 - \exp(-t)]$.

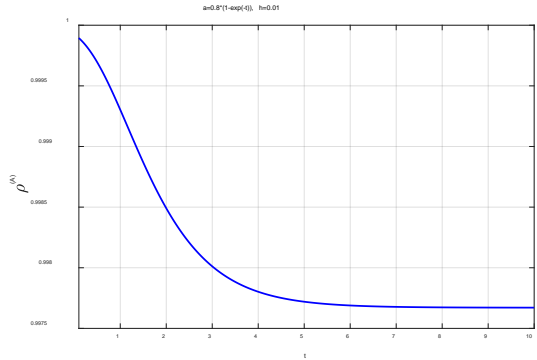

Fig. 7. Spectral radius in Example 2, $a = 0.8[1 - \exp(-t)]$.

## 3.2. Implicit VO-DF equations

We consider first the case where the order of the fractional derivative depends on the function $u(t)$ describing the response of the system. Then Caputo type VO FD of a function $u(t)$ reads

$$D^{\alpha(u)}u(t) = \frac{1}{\Gamma[1-\alpha(u)]}\int_0^t (t-x)^{-\alpha(u)}\dot{u}dx, \qquad 0 \leq \alpha(u) < 1, \ t > 0 \tag{3.12}$$

It is approximated again using Eq. (2.4)

$$(D^{\alpha(t,u)}u)_n = c_1^n(u_n)\dot{u}_1^m + c_2^n(u_n)\dot{u}_2^m \ldots + c_n^n(u_n)\dot{u}_n^m \tag{3.13}$$

where now

$$c_1^n(u_n) = \frac{1}{\Gamma[1-\alpha(u_n)]} \int_0^h (nh-x)^{-\alpha(u_n)} dx$$

$$c_2^n(u_n) = \frac{1}{\Gamma[1-\alpha(u_n)]} \int_h^{2h} (nh-x)^{-\alpha(u_n)} dx \tag{3.14}$$

…

$$c_r^n(u_n) = \frac{1}{\Gamma[1-\alpha(u_n)]} \int_{(r-1)h}^{rh} (nh-x)^{-\alpha(u_n)} dx, \quad r=1,2,\ldots,n \quad n=1,2,\ldots,N$$

Apparently, the coefficients in Eq. (3.13) depend on $u_n$.

Integration of the general term in Eqs. (3.14) gives

$$c_r^n(u_n) = \frac{1}{\Gamma[1-\alpha(u_n)]} \int_{(r-1)h}^{rh} (nh-x)^{-\alpha(u_n)} dx$$
$$= \frac{1}{\Gamma[1-\alpha(u_n)]} \frac{h^{1-\alpha(u_n)}}{\alpha(u_n)-1}[(n-r)^{1-\alpha(u_n)} - [n-r+1]^{1-\alpha(u_n)}] \tag{3.15}$$

$r=1,2,\ldots,n \quad n=1,2,\ldots,N$

Similarly, if the fractional derivative depends on the velocity $\dot{u}(t)$, Eq. (3.12) reads

$$D^{\alpha(\dot{u})}u(t) = \frac{1}{\Gamma[1-\alpha(\dot{u})]} \int_0^t (t-x)^{-\alpha(\dot{u})} \dot{u} dx, \quad 0 \leq \alpha(\dot{u}) < 1, \; t > 0 \tag{3.16}$$

and the counterpart of Eq. (3.15) is

$$c_r^n(\dot{u}_n) = \frac{1}{\Gamma[1-\alpha(\dot{u}_n)]} \frac{h^{1-\alpha(\dot{u}_n)}}{\alpha(\dot{u}_n)-1}[(n-r)^{1-\alpha(\dot{u}_n)} - [n-r+1]^{1-\alpha(\dot{u}_n)}] \tag{3.17}$$

$r=1,2,\ldots,n \quad n=1,2,\ldots,N$

To develop the solution procedure for the implicit VO FDE, we consider again the initial value problem (3.1), in which the VO-FD is defined by Eq. (3.12). Applying this equation for $t=t_n$ and using Eq. (3.13) to approximate the VO-FD, we obtain

$$a_1 q_n + a_2[c_1^n(u_n)\dot{u}_1^m + c_2^n(u_n)\dot{u}_2^m \ldots + c_n^n(u_n)\dot{u}_n^m] + a_3 u_n = p_n \tag{3.18}$$

which is combined with Eqs. (3.2) and (3.3) and solved successively for $n=1,2,$ to yield the solution $u_n$ and the derivatives $\dot{u}_n$, $\ddot{u}_n = q_n$ at instant $t = nh \leq T$. Apparently, Eq. (3.18) is nonlinear. To obtain the solution, it is convenient to combine Eqs. (3.2) and (3.3)

$$\begin{bmatrix} -h & 1 \\ 1 & 0 \end{bmatrix} \begin{Bmatrix} \dot{u}_n \\ u_n \end{Bmatrix} = \begin{bmatrix} 0 & 1 \\ 1 & 0 \end{bmatrix} \begin{Bmatrix} \dot{u}_{n-1} \\ u_{n-1} \end{Bmatrix} + \begin{bmatrix} -\frac{c_1}{2} \\ \frac{c_2}{2} \end{bmatrix} q_n + \begin{bmatrix} -\frac{c_1}{2} \\ \frac{c_2}{2} \end{bmatrix} q_{n-1} \tag{3.19}$$

and solve them for $\dot{u}_n$ and $u_n$. Then substituting in Eq. (3.18) results in a nonlinear algebraic equation for $q_n$. For $n = 1$, the value $q_0$ appears in the right hand side of Eq. (3.19). This value can be obtained from Eq. (3.4) for $t = 0$.

The solution algorithm is presented in Table 1. On the base of this algorithm, a computer program has been written in Matlab for the solution of the initial value problem (3.1), when the. Apparently, this algorithm applies also when VO FD is explicit.

TABLE 1. SOLUTION ALGORITHM

| | |
|---|---|
| A. | Data |
| | Read: $a_1, a_2, a_3, p(t), u_0, \dot{u}_0, t_{tot}$ |
| B. | Initial computations |
| | 1. Choose: $h := \Delta t$ and compute $N = ceiling\ [t_{tot}\ /\ h]$ |
| | 2. Compute $c_1 = h^2/2 \quad c_2 = h \quad q_0 = [p_0 - a_2(D^a u)_0 - a_3 u_0)]\ /\ a_1$ |
| C. | Compute solution |
| | for $n := 1$ to $N$ solve for $\{q_n\ \dot{u}_n\ u_n\}^T$ the system of the nonlinear algebraic equations : |
| | $a_1 q_n + a_2(c_1^n \dot{u}_1^m + c_2^n \dot{u}_2^m \ldots + c_n^n \dot{u}_n^m) + a_3 u_n = p_n$ |
| | $\begin{bmatrix} -h & 1 \\ 1 & 0 \end{bmatrix} \begin{Bmatrix} \dot{u}_n \\ u_n \end{Bmatrix} = \begin{bmatrix} 0 & 1 \\ 1 & 0 \end{bmatrix} \begin{Bmatrix} \dot{u}_{n-1} \\ u_{n-1} \end{Bmatrix} + \begin{bmatrix} -\frac{c_1}{2} \\ \frac{c_2}{2} \end{bmatrix} q_n + \begin{bmatrix} -\frac{c_1}{2} \\ \frac{c_2}{2} \end{bmatrix} q_{n-1}$ |
| | end |

**Example 3**

Solve the initial value problem

$$a_1 \ddot{u} + a_2 D^{a(t)} u + a_3 u = p(t), \quad a(t) = d - k \tanh(|\dot{u}||) \tag{1}$$

Assume: $a_1 = 1$, $a_2 = 2\xi\omega$, $a_3 = \omega^2$, $\xi = 0.1$, $\omega = 2$, $p(t) = 0$, and $u_0, \dot{u}_0, d, k$ as follows:

(i) $d = 0.9999, k = 10e - 10 \quad u_0 = 0,\ \dot{u}_0 = 1$. In this case it is $a(t) \approx 1$, hence as anticipated $D^{a(t)} u \approx \dot{u}$ and the computed solution approximates the exact solution (2) in Example 2. This is shown in Fig 8.

(ii) $d = 1e - 10, k = 1e - 10$, $u_0 = 0,\ \dot{u}_0 = 1$. In this case, it is $a(t) \approx 0$, and the computed solution approximates the exact solution (3) in Example 2. This is shown in Fig 9.

(iii) When $u_0 = 0,\ \dot{u}_0 = 10$, compute the response for: $a(t) = 1; a = 1 - 0.5\tanh(|\dot{u}|)]$. The computed results are plotted in Fig. 10.

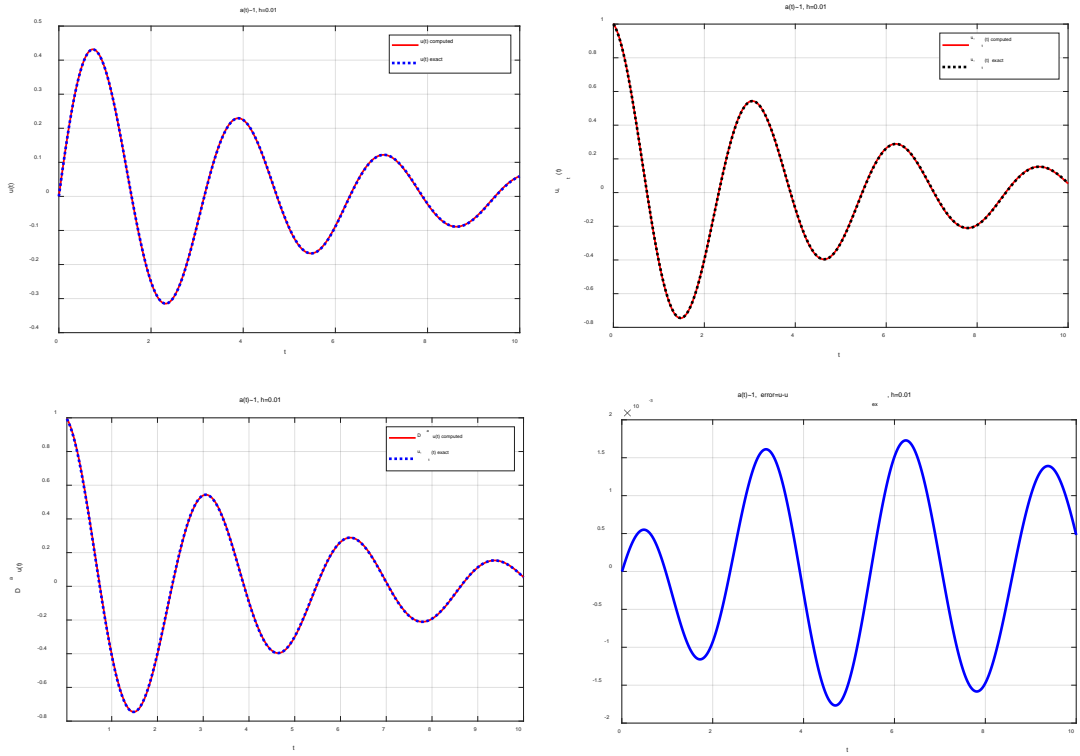

Fig. 8 Results in Example 3: case (i)

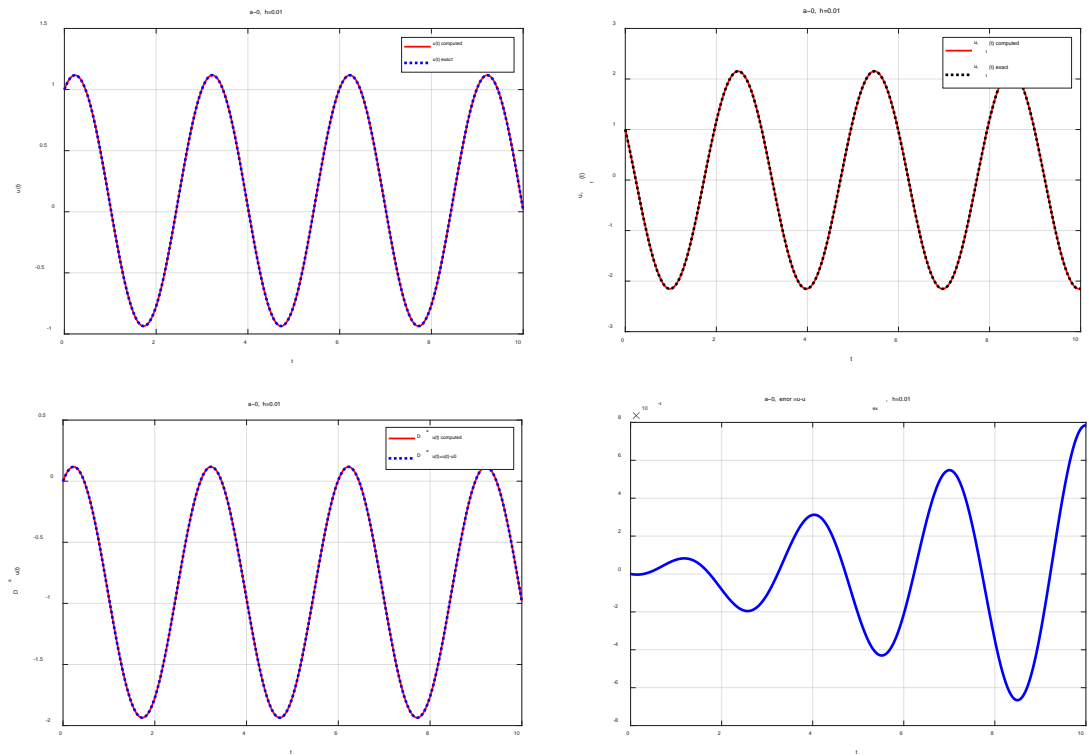

Fig. 9 Results in Example 3: case (ii)

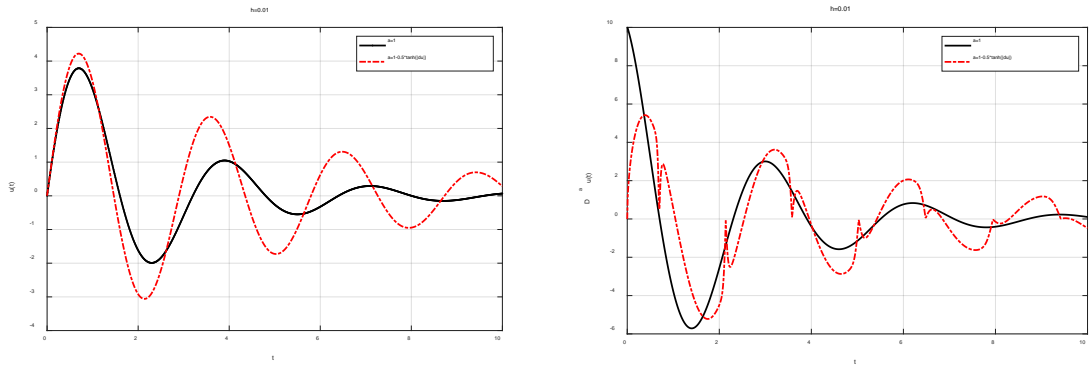

Fig. 10 Results in Example 3: case (iii)

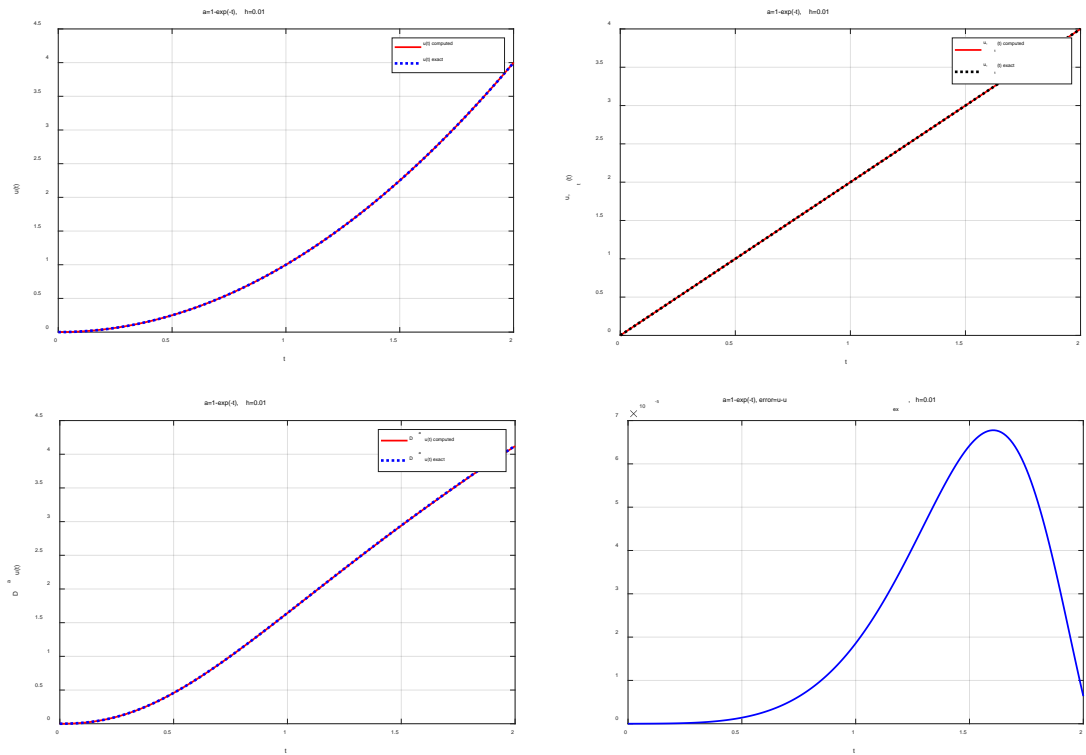

Figure 11. Results in Example 4.

# 4    Nonlinear VO-FD Equations

The solution is obtained using the same algorithm as in linear implicit VO FDEs.

Example 4

**arXiv:1802.00519 [math.NA]**

The numerical scheme is employed to solve the initial value problem for the fractional Duffing oscillator

$$\ddot{u} + 0.2 D^{a(t)} u + u + u^3 = p(t), \quad a = 1 - \exp(-t) \tag{1}$$

$$u_0 = 0, \quad \dot{u}_0 = 0 \tag{2}$$

For $p(t) = 2 + t^2 + t^6 + 0.2 \dfrac{1}{\Gamma(1-a)} \dfrac{2\exp((2t)t^{\exp(-t)+1}}{\exp(t)+1}$, Eq (1) admits an exact solution $u_{ex}(t) = t^2$. The computed results are plotted in Fig. 11.

## 5  VO-FD Equations with variable coefficients

So far we have developed the method for the solution of VO FDEs with constant coefficients. Obviously, if the coefficients $a_1, a_2, a_3$ are functions of the independent variable $t$, the previously described solution procedures remain the same except that the coefficients $a_1, a_2, a_3$ are evaluated in each step. In the following, the effectiveness of the method is demonstrated by solving the initial value problem in the example bellow

**Example 5**

Solve the VO FDE

$$(1+t^2)2 + 0.1 t^{1/2} D^a u + (10 + e^{-t}) u = p(t), \quad a(t) = 1 - 0.5\exp(-t) \tag{1}$$

$$u_0 = 1, \quad \dot{u}_0 = 1 \tag{2}$$

For $p(t) = [(1+t^2) + 0.01 t^{1/2} (\Gamma(1-a,0) - \Gamma(1-a,t))/\Gamma(1-a)) + (10 + e^{-t})] e^t$, Eq (1) admits an exact solution $u_{ex}(t) = e^t$. The computed solution is plotted in Fig. 12 as compared with the exact one.

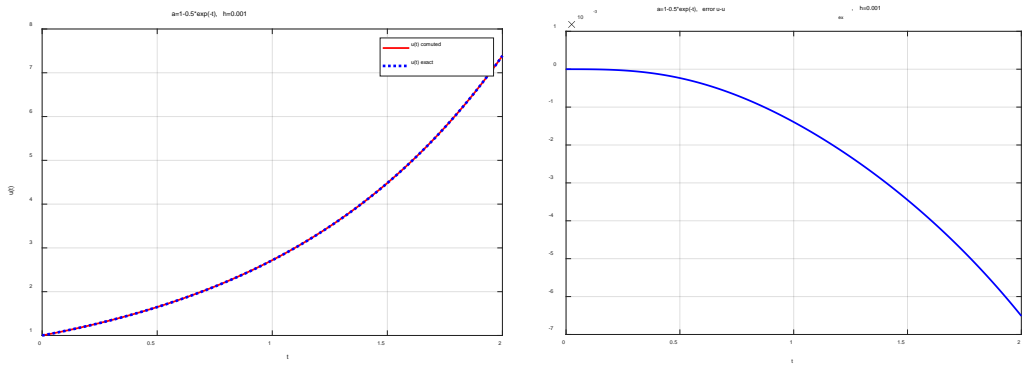

Figure 12. Results in Example 5.

## 5. Conclusions

A numerical method has been developed for solving VO fractional differential equations. The key to the solution procedure is the establishment of a robust simple numerical expression approximating the VO-FD. Then this expression is employed to develop effective numerical schemes for the solution of linear and nonlinear VO-FDEs, both explicit and implicit. The procedure applies also to VO-FDEs with variable coefficients. The method is demonstrated by solving the second order VO FDE, linear and nonlinear, with constant as well as with variable coefficients. However, it can be straightforwardly extended to higher order VO-FDEs and in conjunction with the AEM (Analog Equation Method) to partial VO-FDEs. The developed numerical scheme is stable, accurate and simple to implement. There is no computational complexity in constructing the solution algorithms, while they are easy to program. The convergence and the accuracy of the method are demonstrated through well corroborated numerical examples. The obtained results validate the effectiveness of proposed method.